\documentclass[]{article}
\usepackage{amsmath}
\usepackage{amsfonts}
\usepackage{amssymb}
\usepackage{graphicx}
\usepackage{lineno}
\usepackage{authblk}
\usepackage{verbatim}
\usepackage{mathtools}

\DeclarePairedDelimiter\floor{\lfloor}{\rfloor}
\usepackage{hyperref}
\hypersetup{
	colorlinks=true,
	linkcolor=blue,
	citecolor=violet,
	filecolor=blue,      
	urlcolor=cyan,
}

\usepackage{color}
\usepackage{tikz}
\usetikzlibrary{tikzmark}
\usepackage{float}
\restylefloat{table}
\newtheorem{theorem}{Theorem}[section]

\newtheorem{lemma}[theorem]{Lemma}

\title{\vspace{-1\baselineskip}Strictly Real Fundamental Theorem of Algebra}
\author{Soham Basu}
\affil[]{Corresponding author: sohambasu6817@gmail.com}

\begin{document}

\maketitle

\begin{abstract}
Without resorting to complex numbers \cite{basu} or advanced topological arguments \cite{Pukhlikov, Pushkar}, we show that any real polynomial of degree $>2$ has a real quadratic factor, which is equivalent to the seminal version of the Fundamental Theorem of Algebra (FTA) \cite{Gauss}. Thus it is established that basic real analysis suffices as the minimal platform to prove FTA. 
\end{abstract}



\section{Introduction}
\begin{figure}[h!]
	\centering
	\includegraphics[trim= 300  10   300 55, clip, width=0.65\textwidth]{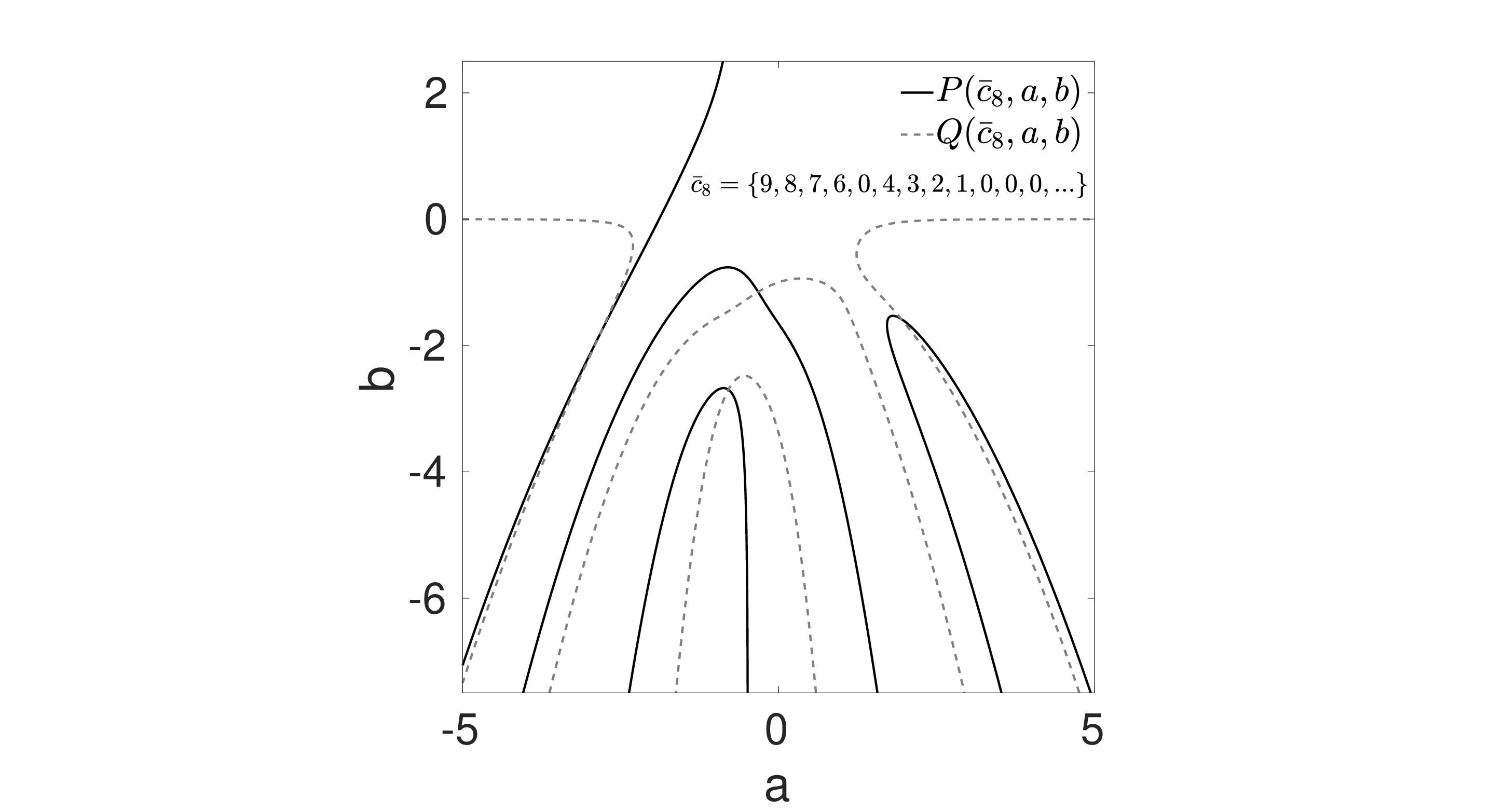}
	\caption{\label{fig:example}$P({\bar{c}}_8,a,b)$ (solid blck lines) and $Q({\bar{c}}_8,a,b)$ (dashed gray lines) corresponding to $f({\bar{c}}_8,~x)=x^8 + 2x^7 + 3x^6 + 4x^4  + 6x^3 + 7x^2 + 8x + 9$, exhibiting interlacing with maximum possible number of (real) roots in the variable $a$ for any fixed $b<-3$. $Q({\bar{c}}_8,a,b=0)$ has no roots in the variable $a$.}
\end{figure}

Let us consider polynomials of the field of real numbers $\mathbb{R}$, which is the minimal field containing the rational number field $\mathbb{Q}$ as a subfield and satisfying the least upper-bound property \cite{Rudin}.
\begin{align}
 f({\bar{c}}_N, ~x) =  \sum_{n=0}^{N} c_{n} x^{n} \in \mathbb{R}
\end{align}
where ${\bar{c}}_N=(c_0 , c_1,c_2,..., c_N, 0,0,0,...) \in \mathbb{R}^{\infty}$ is fixed, with $c_0 \ne 0$ (otherwise divide by $x$ and start redifining ${\bar{c}}_N$ from $c_1$), $c_N \ne 0, ~c_{k}=0~\forall k>N$. $x \in \mathbb{R}$ is independent. Such a description is useful for representing any polynomial as a unique power series. The seminal version of FTA only stated factorizability of any such polynomial into a product of linear and quadratic polynomials \cite{Gauss}, which can be rephrased to an even more compact statement:

\begin{theorem} \label{ch:theorem1}
For $N>2$, given any ${\bar{c}}_N, ~\exists A, B \in \mathbb{R}$ and ${\bar{d}}_{N-2}$ such that 
\begin{align}
f({\bar{c}}_N, ~x) = (x^2-Ax-B) ~ f({\bar{d}}_{N-2}, ~x)
\end{align}
\end{theorem}
We will prove this by looking at the remainder of the polynomial long division of $f(\bar{c}_N, ~x)$ by $(x^2-ax-b)$, where $a,b\in \mathbb{R}$. Given any ${\bar{c}}_N$ and $a,b \in \mathbb{R}$, $\exists$ unique ${\bar{c'}}_{N-2}$, $P({\bar{c}}_N, a,b)\in \mathbb{R}$ and $Q({\bar{c}}_N, a,b)\in \mathbb{R}$ that satisfy
\begin{align} \label{eqn:def}
f({\bar{c}}_N, ~x) = (x^2-ax-b) ~ f({\bar{c'}}_{N-2}, ~x)+ P({\bar{c}}_N, a,b) ~ x + Q({\bar{c}}_N, a,b)
\end{align}
$(x^2-ax-b)=f({\bar{c}}_2, x)$, where ${\bar{c}}_2=(-b, -a, 1, 0,0,0,...)$.~ Both $P(\bar{c}_N, a,b)$ and $Q(\bar{c}_N,a,b)$ are bivariate polynomials in variables $a$ and $b$. By studying their interlacing, we will show the existence of $A,B \in \mathbb{R}$ that satisfy $P({\bar{c}}_N, A,B)=0$ and $Q({\bar{c}}_N, A,B)=0$. The crucial steps will be to show the following:
\begin{enumerate}
	\item $\exists ~ \bar{b}<0$ such that for any fixed $b<\bar{b}$, $P(\bar{c}_N, a,b)$ and $Q(\bar{c}_N, a,b)$ have maximum possible number of (real) roots in the variable $a$ which are interlacing each other. For brevity, this is abbreviated as \textit{interleaving}.
	\item For some $b>\bar{b}$, \textit{interleaving} fails either due to the failure in interlacing of the roots, or due to existence of smaller number of (real) roots of $P(\bar{c}_N, a,b)$ and $Q(\bar{c}_N, a,b)$ than their individual degrees in $a$.
\end{enumerate}
Figure \ref{fig:example} illustrates these conditions for a particular example polynomial.

In order to show such interleaving, we will reduce $P({\bar{c}}_N, a,b)$ and $Q({\bar{c}}_N, a,b)$ to simpler formulae (equation \ref{eqn:simpler}), using the following lemmas: \newline

%\vspace{1\baselineskip}
\begin{lemma}{Linearity}\label{lemma_linearity}
	\begin{align} \label{eqn:linear}
	P({\bar{c}}_N, a,b) = \sum_{n=0}^{N} c_n P_n(a,b) \quad \textrm{and} \quad Q({\bar{c}}_N, a,b) = \sum_{n=0}^{N} c_n Q_n(a,b)
	\end{align}
\end{lemma}
%\textbf{Lemma 1: Linearity} 

%$P({\bar{c}}_N, a,b) = \sum_{n=0}^{N} c_n P_n(a,b)$ and $Q({\bar{c}}_N, a,b) = \sum_{n=0}^{N} c_n Q_n(a,b)$\\
\vspace{1\baselineskip}
where $f(\bar{\epsilon}_n,~x)=x^n$, $P_n(a,b)=P(\bar{\epsilon}_n,a,b)$ and $Q_n(a,b)=Q(\bar{\epsilon}_n,a,b)$.
%, and $\bar{\epsilon}_n=(0,0,...,\epsilon_n=1, 0,0,0,...) \in \mathbb{R}^\infty$. To clarify, $f(\bar{\epsilon}_n,~x)=x^n$, and $\forall n,~  \exists \bar{e}_{n-2}$ such that $x^n = (x^2-ax-b) f(\bar{e}_{n-2}, x) + P_n(a,b) ~x + Q_n(a,b)$. \newline

\textbf{Proof:} Given ${\bar{c}_N}$ and ${\bar{d}}_M$, by definition $\exists {\bar{c'}_{N-2}}$ and $\exists {\bar{d'}}_{M-2}$ such that
\begin{align} \label{eqn:combine}
f({\bar{c}}_N, ~x) &= (x^2-ax-b) ~ f({\bar{c'}}_{N-2}, ~x)+ P({\bar{c}}_N, a,b) ~ x + Q({\bar{c}}_N, a,b) \nonumber \\
f({\bar{d}}_M, ~x) &= (x^2-ax-b) ~ f({\bar{d'}}_{M-2}, ~x)+ P({\bar{d}}_M, a,b) ~ x + Q({\bar{d}}_M, a,b) \nonumber \\
%& \nonumber \\
\implies \alpha_1 f({\bar{c}}_N, ~x) & + \alpha_2 f({\bar{d}}_M, ~x) \nonumber \\
& = (x^2-ax-b) ~ \left[ \alpha_ 1 f({\bar{c'}}_{N-2}, ~x) + \alpha_2 f({\bar{d'}}_{M-2}, ~x) \right] \nonumber \\
& + \left[ \alpha_ 1 P({\bar{c}}_{N}, a,b) + \alpha_2 P({\bar{d}}_{M}, a,b) \right] ~x \nonumber \\
& + \left[ \alpha_ 1 Q({\bar{c}}_{N}, a,b) + \alpha_2 Q({\bar{d}}_{M}, a,b) \right] 
\end{align}
But 
$\alpha_1 f({\bar{c}}_N, ~x) + \alpha_2 f({\bar{d}}_M, ~x)
=\alpha_1 \sum_{n=0}^{N} c_{n} x^{n}+ \alpha_2 \sum_{m=0}^{M} d_{m} x^{m}
=\sum_{k=0}^{max(M,N)} \left[ \alpha_1 c_{k} + \alpha_2 d_{k} \right] x^k = f(\alpha_1 {\bar{c}}_N + \alpha_2 {\bar{d}}_M, ~x)
$. $\exists {\bar{s}}'_{\max\{N,M\}-2}$ such that
\begin{align} \label{eqn:define}
f(\alpha_1 {\bar{c}}_N + \alpha_2 {\bar{d}}_M, ~x) &= (x^2-ax-b) ~ f({\bar{s}}'_{\max\{N,M\}-2}, ~x) \nonumber \\
&+ P(\alpha_1 {\bar{c}}_N + \alpha_2 {\bar{d}}_M, a,b) ~ x + Q(\alpha_1 {\bar{c}}_N + \alpha_2 {\bar{d}}_M, a,b) 
\end{align}
Since polynomial division gives unique remainder, the coefficients of the linear residues in equations \ref{eqn:combine} and \ref{eqn:define} must be equal. Finally the linear representation is written in terms of $P(\bar{\epsilon}_n,a,b)$ and $P(\bar{\epsilon}_n,a,b)$ to obtain equation \ref{eqn:linear}. $\square$ \newline 
%\vspace{\baselineskip}

\begin{lemma}{Recursion} \label{eqn:resursion}
	\begin{align} 
	P_{n+2}(a,b) = a P_{n+1}(a,b) + b P_{n}(a,b) \nonumber \\
	Q_{n+2}(a,b) = a Q_{n+1}(a,b) + b Q_{n}(a,b)
	\end{align}
\end{lemma}

%\textbf{Lemma 2: Recursion} 

\textbf{Proof:} For any $(x^2-ax-b) ~ f(\bar{c}_N) = f(\bar{c}^*_{N+2})$, by the definition of residues in equation \ref{eqn:def}, $P(\bar{c}^*_{N+2},a,b)=0$ and $Q(\bar{c}^*_{N+2},a,b)=0$. Applying lemma \ref{lemma_linearity} on
\begin{align} \label{eqn:chain}
x^{n+2} &= (x^2-ax-b) ~ x^n + a x^{n+1} + b x^n \nonumber \\
&\implies P_{n+2}(a,b)= 0 + a P_{n+1}(a,b) + b P_n(a,b) \nonumber \\
&\text{ and } ~ Q_{n+2}(a,b)= 0 + a Q_{n+1}(a,b) + b Q_n(a,b) \qquad \square
\end{align}

%\vspace{1\baselineskip}
%\textbf{Lemma 3: Interrelation} 
\begin{lemma}{Interrelation}
	\begin{align}\label{eqn:interrelation}
	Q_{n+1}(a,b)=bP_n(a,b)
	\end{align}
\end{lemma}
\textbf{Proof:} By induction. The relation explicitly holds for $n=0,1,2$.

\begin{table}[H]
	\centering
	\begin{tabular}{| c | c | c |}
		\hline
		$n$		& $P_n(a,b)$	 	& $Q_n(a,b)$	 	\\ 
		\hline 
		&				&				\\ [-2ex]
		0		& 0	\tikzmark{a1}& 1				\\
		1		& 1	\tikzmark{b1}& \tikzmark{a2}0				\\
		2		& $a$	\tikzmark{c1}& \tikzmark{b2} $b$				\\
		3		& $a^2+b$ & \tikzmark{c2} $ab$			\\	
		\hline	
	\end{tabular}
\begin{tikzpicture}[overlay, remember picture, yshift=0\baselineskip, shorten >=.5pt, shorten <=.5pt]
\draw [->] ([yshift=.8pt]{pic cs:a2}) -- ({pic cs:a1});
\draw [->] ([yshift=.8pt]{pic cs:b2}) -- ({pic cs:b1});
\draw [->] ([yshift=.8pt]{pic cs:c2}) -- ({pic cs:c1});
\end{tikzpicture}
\end{table}

Assume that equation \ref{eqn:interrelation} holds until some $n>2$. For $n+1$, lemma \ref{eqn:resursion} gives
\begin{align}
Q_{n+2}(a,b)  &= a Q_{n+1}(a,b) + b Q_{n}(a,b) \nonumber \\ 
&= a b P_{n}(a,b) + b^2 P_{n-1}(a,b) \nonumber \\
&= b \left[ a P_{n}(a,b) + b P_{n-1}(a,b) \right] \nonumber \\
&= b P_{n+1}(a,b)  \qquad \square
\end{align}

%\vspace{\baselineskip}

\section{Interleaving using simpler representation}

Using the lemmas, we arrive at a simpler representation of equation \ref{eqn:def}:
\begin{align} \label{eqn:simpler}
f({\bar{c}}_N, ~x) &= (x^2-ax-b) ~ f({\bar{c'}}_{N-2}, ~x) \nonumber \\
&+ x~ \sum_{n=1}^{N} c_n P_n(a,b) + b~ \sum_{n=1}^{N} c_n P_{n-1}(a,b) + c_0
\end{align}
Therefore we only need to look at the interleaving of $\sum_{n=1}^{N} c_n P_n(a,b)$ and $b~ \sum_{n=1}^{N} c_n P_{n-1}(a,b) + c_0$. Since the extra variable $b$ and constant $c_0$ renders direct comparison of their roots intractable, we resort to comparing roots of the following polynomials similar in structure to the famous Sturm chains \cite{Sturm}:
\begin{align} \label{eqn:Sturm}
h_{m}(a,b)=\sum_{n=m}^{N} c_{n} P_{n-m}(a,b) 
\end{align}
$h_{m}(a,b)$ inherits the following recursion formula from lemma 2:
\begin{equation} \label{eqn:recursion}
h_{m}(a,b) = a h_{m+1}(a,b) + b h_{m+2}(a,b) + c_{m+1} 
\end{equation}
$\sum_{n=1}^{N} c_n P_n(a,b) = h_0(a,b)$ and $b~ \sum_{n=1}^{N} c_n P_{n-1}(a,b) + c_0 = b h_1(a,b) + c_0$.

 %\clearpage
 \vspace{1\baselineskip}
 %\textbf{Theorem 2:}
 \begin{theorem} \label{ch:interleaving_theorem} 
 	For any given ${\bar{c}}_N, ~\exists b_{{\bar{c}}_N}<0$  such that for any fixed $b<b_{{\bar{c}}_N}$ , and for all $0\le m \le N-2$, $\big( h_{m}(a,b),~ h_{m+1}(a,b) \big)$ are interleaving. Recapitulating for the sake of clarity, for any fixed $b<b_{{\bar{c}}_N}$,
 	\begin{itemize}
 		\item $h_{m}(a,b)$ has $N-m-1$ distinct (real) roots in the variable $a$, say $\alpha_1^{(m)}(b)< \alpha_2^{(m)}(b)<...<\alpha_{N-m-1}^{(m)}(b)$ in ascending order. 
 		\item All pairs $\big( h_{m}(a,b),~ h_{m+1}(a,b) \big)$ have interlacing roots $\implies \alpha_1^{(m)}(b)<  \alpha_1^{(m+1)}(b) < \alpha_2^{(m)}(b) < \alpha_2^{(m+1)}(b) <...<\alpha_{N-m-2}^{(m+1)}(b) <\alpha_{N-m-1}^{(m)}(b)$.
 	\end{itemize}  
 \end{theorem} 
 %For brevity, let us call these conditions together to be \textit{interleaving}.

 %Fisk\cite{Fisk} showed that is $h_{m}$ were general polynomials of degree $N-m-1$, the recursive equation in (6) without $c_{m}$ enforces existence of maximum number of real roots and interlacing of the roots for each pair $(h_{m+1}(a,b),h_{m}(a,b)$, provided that $(h_{N-3}(a,b),h_{N-2}(a,b))$ interlace roots. 
 \vspace{\baselineskip}
 \textbf{Proof:} By induction and contradiction. We can explicitly verify for $m=N-3$ and $m=N-2$
 \begin{align} \label{eqn:explicit}
 h_{N-3}(a,b) &= c_N P_3(a,b)+c_{N-1}P_2(a,b) + c_{N-2}P_1(a,b) \nonumber \\
 &= c_N (a^{2} + b) +c_{N-1}a+c_{N-2} \nonumber \\
 h_{N-2}(a,b) &= c_N P_2(a,b) +c_{N-1} P_1(a,b) =c_N  a + c_{N-1}
 \end{align}
 
 $h_{N-2}(a,b)=0$  only at $\alpha^{(N-2)}_1(b) = - \frac{c_{N-1}}{c_N}$. For any fixed $b< b_{N-3}=- \frac{c_{N-2}}{c_N}$, $h_{N-3}(a,b$) has distinct (real) roots. When $c_{N-1} \ne 0$, these roots are $-\frac{c_{N-1}}{2 c_N} (1 \pm \sqrt{1 - D} )$, where $D={\frac{4 c_N c_{N-2} }{ c_{N-1}^2 } +  \frac{4 c_N^2 }{ c_{N-1}^2 } b}$. When $c_{N-1} = 0$, the roots are $\pm \sqrt{-\frac{c_{N-2}}{c_N}-b}$. 
 %Designating the roots of $h_{N-3}(a,b)$ as $\{ \alpha^{(N-3)}_l(b),~l=1,2 \}$ in ascending order (For example 
Either $\alpha^{(N-3)}_1(b)=-\frac{c_{N-1}}{2 c_N} (1 + \sqrt{1 - D} )$ when $\frac{c_{N-1}}{c_N}>0$, or $\alpha^{(N-3)}_1(b)=-\frac{c_{N-1}}{2 c_N} (1 - \sqrt{1 - D} )$ when $\frac{c_{N-1}}{c_N}<0$. By directly comparing the formulae of roots, we get interleaving of $\big( h_{N-3}(a,b), ~h_{N-2}(a,b) \big)~ \forall b<b_{N-3}$ :
 \begin{align}
 \alpha^{(N-3)}_1(b) <  \alpha^{(N-2)}_1(b) < \alpha^{(N-3)}_2(b) \nonumber 
 \end{align} 
 \vspace{-\baselineskip} 
 \begin{figure}[h!]
 	\centering
 	\includegraphics[trim= 262  10   268 62, clip, width=0.99\textwidth]{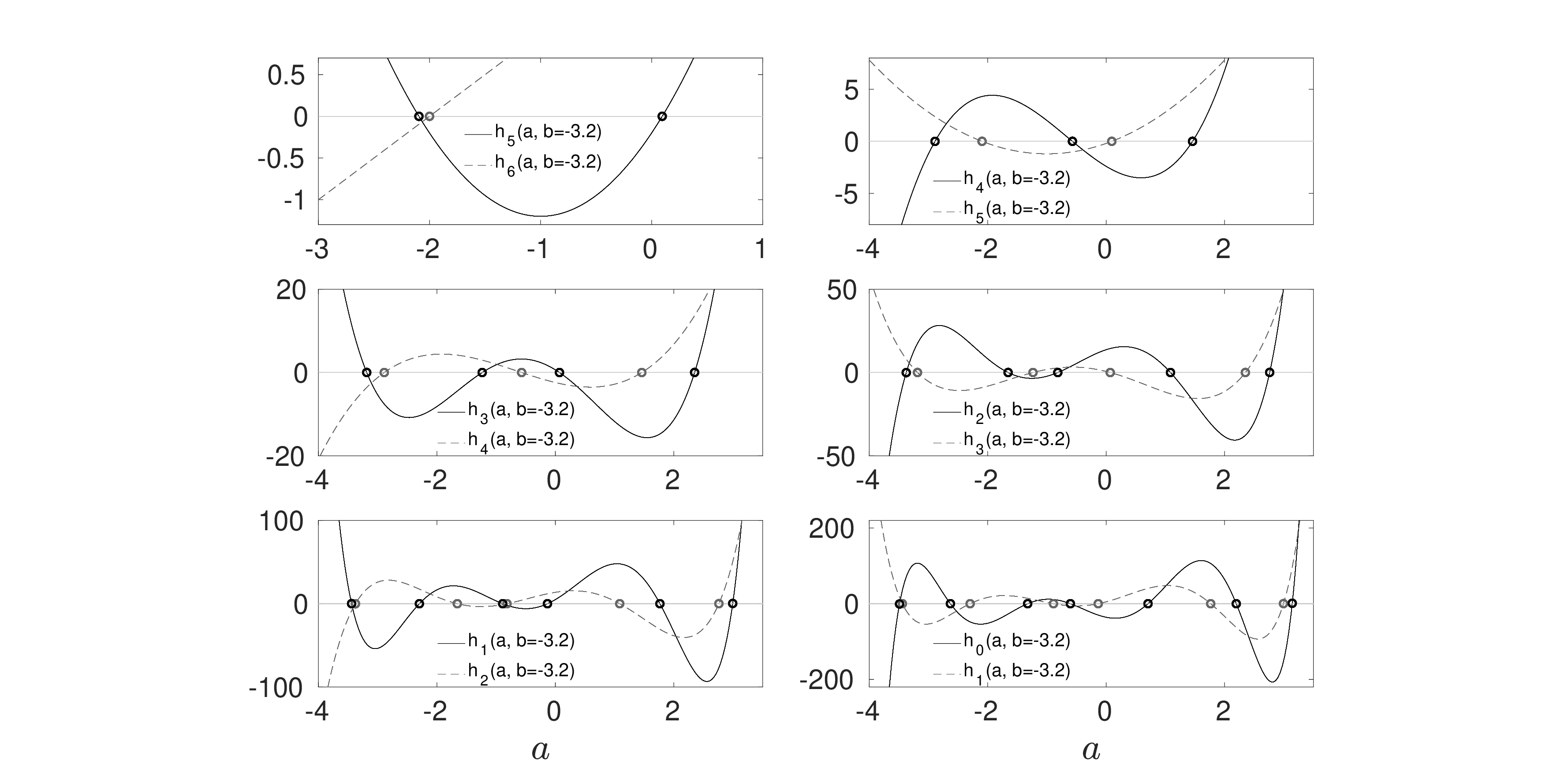}
 	\caption{\label{fig:interleaving} For ${\bar{c}}_8=(9, 8, 7, 6, 0 , 4, 3, 2, 1, 0,0,0,...)$, computed example of interleaving of $\big( h_m(a, b=-3.2), ~ h_{m+1}(a,b=-3.2) \big)$ for $0 \le m \le N-3$. Black circles illustrate $\alpha_l^{(m)}(b=-3.2)$ and gray circles illustrate $\alpha_l^{(m+1)}(b=-3.2)$. }
 \end{figure}

%$h_{N-3}(a,b$) starts to have distinct (real) roots when $b< b_{N-3}=( \frac{c_{N-1}^{2}}{4 c_N^2} - \frac{c_{N-2}}{c_N})$, $-\frac{c_{N-1}}{2 c_N} - \frac{c_{N-1}}{2 c_N} \sqrt{1- \frac{4 c_N c_{N-2} }{ c_{N-1}^2 } -  \frac{4 c_N^2 }{ c_{N-1}^2 } b}$ and $-\frac{c_{N-1}}{2 c_N} + \frac{c_{N-1}}{2 c_N} \sqrt{1- \frac{4 c_N c_{N-2} }{ c_{N-1}^2 } -  \frac{4 c_N^2 }{ c_{N-1}^2 } b}$

%both $h_{N-3}(a,b)$ and $h_{N-2}(a,b)$ have real roots which interlace. So $b_{{\bar{c}}_N}< b_{N-3}$ suffices for $\Big( h_{N-3}(a,b),~ h_{N-2}(a,b) \Big)$.

%\vspace{\baselineskip}
Now let us assume that for the given ${\bar{c}}_N, ~ \exists b_{m+1}<0$ such that for any fixed $b<b_{m+1}$, the pairs $\big( h_{k}(a,b),~ h_{k+1}(a,b) \big)$ interleave for $m+1 \le k \le N-3$. %For $b<b_{m+1}$, let us designate by $\{ \alpha^{(m+1)}_{l}(b): 1 \le l \le N-m-2 \}$ the roots of $h_{m+1}(a,b)$ in ascending order: $\alpha^{(m+1)}_{1}(b)<\alpha^{(m+1)}_{2}(b))<...<\alpha^{(m+1)}_{N-m-2}(b)$. 
Interleaving implies that $h_{m+2}\big(\alpha^{(m+1)}_{l}(b),b \big) \quad h_{m+2} \big( \alpha^{(m+1)}_{l+1}(b),b \big) <0, ~ \forall l$. 
\vspace{\baselineskip}

In the particular case of $c_{m+1}=0$, equation \ref{eqn:recursion} implies that for any such $b$, 
\begin{align}
& h_{m}(\alpha^{(m+1)}_{l}(b),b) \quad h_{m}(\alpha^{(m+1)}_{l+1}(b),b) \nonumber \\
= & \Big( 0 + b h_{m+2}(\alpha^{(m+1)}_{l}(b),b) \Big) ~ \Big( 0 + b h_{m+2}(\alpha^{(m+1)}_{l+1}(b),b) \Big) <0
\end{align}
Applying intermediate value theorem establishes the existence of $\alpha^{(m)}_{l}(b) \in  \big(\alpha^{(m+1)}_{l-1}(b), ~ \alpha^{(m+1)}_{l}(b) \big), ~ \forall~ 2 \le l \le N-m-2$ such that $h_{m}(\alpha^{(m)}_{l}(b),b)=0$.
Since $h_{m}(a,b)$ and $h_{m+2}(a,b)$ both have the same leading coefficient $c_N$ and differ by degree 2 in the highest power of $a$, they take same signs at $\pm \infty$: $\lim_{a \rightarrow \pm \infty} h_{m}(a,b) ~ h_{m+2}(a,b) = \infty$. In addition, due to the negative sign of $b<b_{m+1}<0$, we have at the extremal roots of $h_{m+1}(a,b)$:
\begin{align}
	& h_{m} \big( \alpha^{(m+1)}_{1}(b),b \big) \quad h_{m+2} \big( \alpha^{(m+1)}_{1}(b),b \big) = b \Big( h_{m+2} \big( \alpha^{(m+1)}_{1}(b),b \big) \Big)^2  <0 \nonumber \\
	& h_{m} \big( \alpha^{(m+1)}_{N-m-2}(b),b \big) \quad h_{m+2} \big( \alpha^{(m+1)}_{N-m-2}(b),b \big) = b \Big( h_{m+2} \big( \alpha^{(m+1)}_{N-m-2}(b),b \big) \Big)^2 <0 \nonumber	
\end{align} 
Due to the location of all its roots in $I(b<b_{m+1})=\big( \alpha^{(m+1)}_{1}(b), \alpha^{(m+1)}_{N-m-2}(b)\big)$, $h_{m+2}(a,b)$ cannot change sign outside $I(b)$. Therefore $h_m(a,b)$ must change sign outside $I(b)$. By the intermediate value theorem (real) roots must exist in both $\big( -\infty, \alpha^{(m+1)}_{1}(b) \big)$ and $\big( \alpha^{(m+1)}_{N-m-2}(b), \infty \big)$. Thus in the particular case of $c_{m+1}=0$ interleaving of the pair $\big( h_{m}(a,b),~ h_{m+1}(a,b) \big)$ is directly established. 
\vspace{\baselineskip}

For the general case of $c_{m+1} \ne 0$, intermediate value theorem is not applicable using equation \ref{eqn:recursion} without $b h_{m+2}(\alpha^{(m+1)}_{l}(b),b)$ somehow dominating $c_{m+1}$. Let us now derive such a sufficient growth property for large negative $b$:

\begin{lemma} \label{ch:theorem_growth} {Growth property} $\qquad \mathrm{For ~ any ~} 1 \le l \le N-m-2$
	\begin{align} \label{eqn:growth}
	\lim_{b \rightarrow -\infty}    \vert h_{m+2}(\alpha^{(m+1)}_{l}(b), b) \vert \rightarrow \infty
	\end{align}
\end{lemma}

\textbf{Proof:} By contradiction. Let us fix an $1 \le l \le N-m-2$. Assuming that lemma \ref{ch:theorem_growth} is false (equivalent to  $\liminf_{b \rightarrow -\infty}    \vert h_{m+2}(\alpha^{(m+1)}_{ l}(b)) \vert < \infty$) implies that $\exists \epsilon >0$ such that we can find at least one unbounded sequence $\left \{ b_{i}: i \in \mathbb{N}, b_{i+1}<b_{i}<0 \right \}$  that satisfies $\vert h_{m+2}(\alpha^{(m+1)}_{l}(b_{i}), b_{i}) \vert < \epsilon$. Consequently from equation \ref{eqn:recursion}
 \begin{align} \label{eqn:bounds}
  h_{m+1}&( \alpha^{(m+1)}_{l}(b_{i}), b_i ) = \alpha^{(m+1)}_{l}(b_{i}) ~ h_{m+2}(\alpha_{m+1}(b_{i}), b_i)\nonumber \\
  & \qquad \qquad \qquad \qquad +  b_i h_{m+3}(\alpha^{(m+1)}_{l}(b_{i}), b_i) + c_{m+2}  \nonumber \\
  \implies& 0 =  h_{m+3}(\alpha^{(m+1)}_{l}(b_{i}), b_{i})+ \frac{\alpha^{(m+1)}_{l}(b_{i})}{b_{i}} ~ h_{m+2}(\alpha_{m+1}(b_{i}), b_i)  
  + \frac{c_{m+2}}{b_{i}} \nonumber \\
  \implies  &
  \vert h_{m+3}(\alpha^{(m+1)}_{l}(b_{i}), b_{i}) \vert \le 
  \vert \frac{\alpha^{(m+1)}_{ l}(b_{i})}{b_{i}} \vert \vert h_{m+2}(\alpha_{m+1}(b_{i}), b_i) \vert 
  +\vert\frac{c_{m+2}}{b_{i}}\vert \nonumber \\
  \implies  &
  \vert h_{m+3}(\alpha^{(m+1)}_{l}(b_{i}), b_{i}) \vert \le 
  \vert \frac{\alpha^{(m+1)}_{ l}(b_{i})}{b_{i}} \vert ~ \epsilon
  +\vert\frac{c_{m+2}}{b_{1}}\vert
 \end{align}

 Writing $h_{m+1}(a,b)$ as a polynomial in $a$,  $h_{m+1}(a,b)=\sum_{k=0}^{N'}g_{N'-k}(b)a^{k}$ with $N'=N-m-2$, and then applying triangle inequality we get $\forall a \ne 0$,
 
\begin{align}
	& \vert h_{m+1}(a,b) \vert \ge  \vert c_N  \vert \vert a  \vert^{N-m-2} \bigg( 1- \sum_{k=1}^{N-m-2}\frac{\vert g_{k}(b) \vert } {\vert c_N  \vert \vert a \vert ^{k}} \bigg) \nonumber \\
	& \implies \vert h_{m+1}(a,b) \vert >0 \quad \forall \vert a \vert > \max_{1 \le k \le N-m-2} \left[ (\frac{ \vert g_{k}(b) \vert}{N \vert c_N  \vert}) ^{\frac{1}{k}} \right] \nonumber \\
	& \implies \vert \alpha^{(m+1)}_{l} (b_{i}) \vert < \max_{1 \le k \le N-m-2} \left[ (\frac{ \vert g_{k}(b) \vert}{N \vert c_N  \vert}) ^{\frac{1}{k}} \right]	\nonumber \\
	& \implies \vert \frac{\alpha^{(m+1)}_{l}(b_{i}) }{b_{i}} \vert < \max_{1 \le k \le N-m-2} \left[ (\frac{ \vert g_{k}(b_{i}) \vert}{N \vert c_N  \vert \vert b_{i} \vert^{k}}) ^{ \frac{1}{k}} \right ]	
\end{align}  
 
 %Thus we get the following bound on $\vert \frac{\alpha_{m} (b_{i})}{b_{i}} \vert$

 Since the appearance of one $b$ while deriving the formula for $h_{m}(a,b)$ from $\sum_{n=0}^{N-m}c_{N-m-n} x^{n}$ is equivalent to continually replacing $x^2$ with $ax+b$ until only linear residue is left is $x$, it can be seen that appearance of one $b$ in such replacement is equivalent to losing $x^2$. This gives $\floor*{k/2}$ as the upper limit of the power of $b$ in $g_{k}(b)$, where $\floor*{y}$ gives the largest integer $\le y \in \mathbb{R}$. Thus 
 \begin{align}
 & \lim_{b_{i} \rightarrow -\infty} \quad \max_{1 \le k \le N-m-2} \left[  \bigg( \frac{ \vert g_{k}(b_{i} ) \vert}{N \vert c_N \vert \vert b_{i} \vert^{k}}  \bigg) ^{ \frac{1}{k}} \right ] \nonumber \\
 & \le 
 \lim_{b_{i} \rightarrow -\infty} \quad \max_{1 \le k \le N-m-2} K \left[  \bigg( \frac{ \vert b_{i} \vert^{ \floor*{k/2}}}{N \vert c_N \vert \vert b_{i} \vert^{k}}  \bigg) ^{ \frac{1}{k}} \right ] \text{ for some } K>0 \nonumber \\
 & \le \lim_{b_{i} \rightarrow -\infty} \frac{K}{(N \vert c_N \vert) ^{\frac{1}{N}}} ~ \vert b_{i} \vert^{\frac{1}{2}}   = 0   \implies \lim_{b_{i} \rightarrow -\infty}  \vert \frac{\alpha^{(m+1)}_{l}(b_{i}) }{b_{i}} \vert = 0
 \end{align} 
 %For our purposes, only the boundedness of $\vert \frac{\alpha^{(m+1)}_{l}(b_{i}) }{b_{i}} \vert <K_m, ~\forall b_i$ suffices. Applied to equation \ref{eqn:bounds}, such boundedness implies that
 Thus we get a strong limiting behavior over $\left \{ b_{i}: i \in \mathbb{N}, b_{i+1}<b_{i}<0 \right \}$:
 \begin{align} \label{eqn:limit}
	 \lim_{b_i \rightarrow -\infty} \vert h_{m+3}(\alpha^{(m+1)}_{l}(b_{i}), b_{i}) \vert \le 
	 \lim_{b_{i} \rightarrow -\infty} \left[ \epsilon   \vert \frac{\alpha^{(m+1)}_{l}(b_{i}) }{b_{i}} \vert
	 +\vert\frac{c_{m+2}}{b_{i}}\vert \right] = 0
 \end{align}
 
Applying limit \ref{eqn:limit} to equation \ref{eqn:recursion} for $m+2$, we get similar limiting behavior
  \begin{align}
 & \lim_{b_i \rightarrow -\infty} \vert h_{m+4}(\alpha^{(m+1)}_{l}(b_{i}), b_{i}) \vert \le  
 \lim_{b_i \rightarrow -\infty} \vert \frac{\alpha^{(m+1)}_{ l}(b_{i})}{b_{i}} \vert
 ~ \lim_{b_i \rightarrow -\infty} \vert h_{m+3} \big (\alpha_l^{(m+1)}(b_{i}), b_{i} \big) \vert  \nonumber \\
 & \qquad \qquad \qquad \qquad +  \lim_{b_i \rightarrow -\infty}\vert \frac{h_{m+2} \big( \alpha_l^{(m+1)}(b_{i}), b_{i} \big) }{b_{i}} \vert
 + \lim_{b_i \rightarrow -\infty} \vert\frac{c_{m+3}}{b_{i}}\vert = 0 
 \end{align}

Continuing sequentially for larger $m$, we deduce for $m=N-3$ and $m=N-2$: 
\begin{align}
	& \lim_{b_i \rightarrow -\infty} h_{N-3}(\alpha^{(m+1)}_{l}(b_{i}), b_{i}) =0  \label{eqn:limit_N-3} \\ 
	& \lim_{b_i \rightarrow -\infty} h_{N-2}(\alpha^{(m+1)}_{l}(b_{i}), b_{i}) =0 \implies \lim_{b_i \rightarrow -\infty} \alpha^{(m+1)}_{l}(b_{i}) = -\frac{c_{N-1}}{c_N} \label{eqn:limit_N-2}
\end{align}

If equations \ref{eqn:limit_N-3}  and and \ref{eqn:limit_N-2} hold true simultaneously, then by equations \ref{eqn:explicit} 
\begin{align} \label{eqn:lims}
& \lim_{b_i \rightarrow -\infty} h_{N-3}(\alpha^{(m+1)}_{l}(b_{i}), b_{i}) = c_N (\frac{c_{N-1}}{c_N})^2 + \lim_{b_i \rightarrow -\infty} c_N b_i - c_{N-1} \frac{c_{N-1}}{c_N} + c_{N-2} \nonumber \\
%= & \lim_{b_i \rightarrow -\infty} h_{N-2}(\alpha^{(m+1)}_{l}(b_{i}), b_{i}) \quad \lim_{b_i \rightarrow -\infty} \frac{ h_{N-2}(\alpha^{(m+1)}_{l}(b_{i}), b_{i}) -c_{N-1} }{c_N}  \nonumber \\
& = \lim_{b_i \rightarrow -\infty} c_N b_i + c_{N-2}= - \infty \frac{c_N}{\vert c_N \vert}
\end{align}
To avoid contradiction between equations \ref{eqn:limit_N-3} and \ref{eqn:lims}, lemma \ref{ch:theorem_growth} must hold. $\square$
Let us define $b_{m+1}<0$ such that $\forall b<b_{m+1}$, $\vert b h_{m+2}(\alpha^{(m+1)}_{l}(b), b) \vert > \vert c_{m+1}\vert$ holds, thereby making intermediate value theorem applicable. Continuing this procedure sequentially over $m$, we get the set $\{ b_{m}: N-3 \le m \le 0 \}$.  $b_{{\bar{c}}_N} = \displaystyle \min_{N-3 \le m \le 0} \{ b_{m} \}$ proves theorem \ref{ch:interleaving_theorem}. $\square$

\vspace{\baselineskip}

For $m=0$, lemma \ref{ch:theorem_growth} implies that $\lim_{b_i \rightarrow -\infty} \vert h_{1}(\alpha^{(0)}_{l}(b_{i}), b_{i}) \vert = \infty$. Thus $\exists ~ \bar{b}<0$ such that $\vert  h_{1}(\alpha^{(0)}_{l}(b), b) \vert > \vert \frac{c_0}{b} \vert$ holds $\forall b< \bar{b}$. Applying the intermediate value theorem proves interleaving of $\big( h_{0}(a,b), ~ h_{1}(a,b) + \frac{c_0}{b} \big)$ for any fixed $b<\bar{b}$. Since for fixed $b$ the roots of $h_{1}(a,b) + \frac{c_0}{b}$ and $b h_{1}(a,b) + {c_0}$ are the same, $\big( h_{0}(a,b), ~ b h_{1}(a,b) + {c_0} \big)$ also interleave for any fixed $b<\bar{b}$. 
% Since $b$ is fixed, this shows that $\big( h_{0}(a,b), bh_{1}(a,b) + c_0 \big)$ interleave roots $\forall b<\bar{b}$.

%Using $\lim_{b \rightarrow -\infty}    \vert h_{1} \big( \alpha^{(0)}_{l}(b), b \big) \vert \rightarrow \infty$ and interleaving of  $(h_{0}(a,b),~ h_{1}(a,b)$, it is straightforward to show existence of some $b_{{\bar{c}}_N} \le \min_{0 \le N-2} b_m$ such that $(h_{0}(a,b), ~  b h_{1}(a,b)+ {c_0})$ interleaves $\forall ~ b < b_{{\bar{c}}_N}$. The values of $b$ which can allow non-interleaving of $(h_{0}(a,b), ~  b h_{1}(a,b)+ {c_0})$ is thus bounded below by $b_{{\bar{c}}_N}$. $b=0$ is explicitly a non-interleaving value, since $b h_{1}(a,b)+ {c_0}$ has no root at $b=0$. By the least upper bound property of real numbers, the non-interleaving set of $(h_{0}(a,b), ~  b h_{1}(a,b)+ {c_0})$ over $b$ has a smallest lower bound, say $\bar{b}<0$. 

%have $N-1$ and $N-2$ interlacing real roots in an unbounded set of negative $b$. This reduces the problem statement to \textit{our method}. \newline
 
 \section{Existence of $A,B$ satisfying $P({\bar{c}}_N,A,B)=0$ and $Q({\bar{c}}_N,A,B)=0$ from interleaving}
Let us consider values of $b$ for which $\big( h_{0}(a,b), ~ b h_{1}(a,b) + {c_0} \big)$~ does not interleave. Clearly this non-interleaving set is bounded below by $\bar{b}$. For small $\vert b \vert $, $h_{1}(a,b)=-\frac{c_0}{b}$ can only be attained either at two values of $a$ or none, due to the boundedness of the coefficients of $h_{1}(a,b)$. Particularly at $b=0$, $bh_{1}(a,b) + c_0=c_0$ is only a constant and therefore has no roots. Thus the non-interleaving set contains some interval around $b=0$. Since the non-interleaving set is non-empty and bounded below, by the least upper-bound property \cite{Rudin} it has an infimum $B \in  \mathbb{R}$ satisfying $\bar{b} \le B <0$. Since the coefficients of both $h_{0}(a,b)$ and $b h_{1}(a,b) + {c_0}$, when written as polynomials in $a$, are bounded for $\bar{b} \le b \le B<0$, the roots of both $h_{0}(a,b)$ and $b h_{1}(a,b) + {c_0}$ in variable $a$ are bounded for these values of $b$. Thus by Bolzano-Weierstrass theorem, roots can be defined at $b=B$ as subsequential limits of roots taken over any sequence $\{ b_k: k \in \mathbb{N}, ~\bar{b} \le b_k < B, ~ b_k < b_{k+1}, ~ \lim_{k \rightarrow \infty} b_k =B \}$. Writing the roots of $b h_{1}(a,b) + {c_0}$ as $\beta_l(b)~ \forall b < B$, we can define limits over any converging subsequences of $ \beta_{l}(b_{k})$ and $\alpha_{l}^{(1)}(b_{k})$ indexed by $\sigma_k$ and $\sigma'_k$.
\begin{align}
\beta_j=\lim_{j \rightarrow \infty} \beta_{l}(b_{\sigma_k}),& \quad 1 \le i \le N-2 \nonumber \\
\alpha_i=\lim_{j \rightarrow \infty} \alpha_{l}^{(1)}(b_{\sigma'_k}),& \quad 1 \le i \le N-1 \nonumber \\
\implies \alpha_1 \le \beta_1 \le \alpha_2 \le ... \le \beta_{N-2} & \le \alpha_{N-1} \qquad  \forall b \le B \label{eqn:atB}
\end{align}

If all the inequalities are strict, then $B$ belongs to the interleaving set. Then $\exists t>0$ such that $(\alpha_i-t, \alpha_i+t)$ and $(\beta_j-t, \beta_j+t)$ are disjoint. This implies 
\begin{align}
& h_0(\alpha_i-t,B) h(\alpha_i+t,B) <0 \nonumber \\
& \big( {B} h_1(\beta_j-t,B)+{c_0} \big)  \big( {B} h_1(\beta_j+t,B)+{c_0} \big)  <0
\end{align} 
Due to continuity of bivariate polynomials, $\exists  ~ 0 <\delta <\frac{-B}{2}$ such that $h_0(a,b)$ and $b h_{1}(a,b) + c_0$ do not change sign in inside circles of radius $\delta$ centered around $(\alpha_i \pm t,B)$ and $(\beta_j \pm t,B)$, respectively. Then $\forall ~ B \le b \le B+{\delta}$
\begin{align}
	&h_0(\alpha_i-t,b) ~ h_0(\alpha_i+t,b) <0 \quad \implies \forall B \le b \le B+{\delta}    \nonumber \\
	&\qquad h_0(a,b) \text{ has at least one root in } (\alpha_i- t, \alpha_i+t) \nonumber \\
	&\big( b h_1(\beta_j- t,b)+{c_0} \big) ~ \big( b h_1(\beta_j+t,b)+{c_0} \big) <0  \quad \implies \forall B \le b \le B+{\delta} \nonumber \\
	&\qquad b h_1(a,b)+{c_0}\text{ has at least one root in } (\beta_j- t, \beta_j+t)  
\end{align} 
But then interleaving holds $\forall ~ B \le b \le B+{\delta}$, which contradicts the definition of $B$ as the infimum of the  non-interleaving set. Therefor at least two neighboring quantities in equation \ref{eqn:atB} must be equal. The value of $A$ is given by these equal quantities at $b=B$, thereby proving theorem \ref{ch:theorem1}. $\blacksquare$
\vspace{\baselineskip}

Additionally, we have proven that we can always find $B<0$ satisfying theorem \ref{ch:theorem1} for any ${\bar{c}}_N, ~N>2$. To obtain complete factorization of $f({\bar{c}}_N, ~x)$, factorization can be continued onward from $f({\bar{d}}_{N-2}, ~x)$ until it halts due to lack in powers of $x$.

\end{document}